\documentclass[a4paper,twoside,final]{article}
\usepackage{amssymb}
\usepackage[latin1]{inputenc}
\usepackage[T1]{fontenc}
\usepackage[final]{lpretex}
\usepackage[all]{xy}
\usepackage{amscd}
\usepackage{title}

\def\g{{\frak g}}
\def\h{\frak h}
\def\b{\frak b}
\def\n{\frak n}
\def\Pow#1{\bigl[[#1]\bigr]}
\newcommand\oA{\bar{A}}
\let\Emph\textit

\begin{document}

\title{Kac-Moody algebras and the cde-triangle}
\author{Torsten Ekedahl}
\address{Department of Mathematics\\
 Stockholm University\\
 SE-106 91  Stockholm\\
Sweden}
\email{teke@math.su.se}

\begin{abstract}
This is the written version of a talk at the conference on ``Non-commutative
geometry and representation theory in mathematical physics'' held in Karlstad,
Sweden, 5--10 July, 2004. In it we show that the duality formula of Rocha-Caridi
and Wallach is a simple consequence of the so called cde-triangle of modular
representation theory. It tries to reflect the attempt of the talk to cater to
the differing backgrounds of its listeners.
\end{abstract}
\maketitle

The purpose of this note is to give a short proof of the so called duality
formula \cite{rocha-caridi82::projec+lie} occurring in the representation theory
of Kac-Moody Lie algebras. The proof uses techniques from the theory of modular
representations, in fact the duality formula turns out to be a special case of
the fact that ``e is the transpose of d'' \cite[III,\S15]{serre78::repres}. More
precisely, the idea is that for a generic weight, the Verma module is
irreducible and one has a semi-simple category of representations. The specific
characters are then obtained from a generic character by specialisation, i.e.,
one constructs Verma modules depending on a parameter such that for the general
value the Verma module is irreducible and for a special value we are in the
particular situation that we are interested in. In this way the modular set up
is over a power series ring in one variable over the base field rather than over
a $p$-local ring.  There is a technical difficulty (which is not serious) in
that the modules over the Kac-Moody algebra that are considered are not finitely
generated over the base ring.
\begin{section}{The cde-triangle}
\label{sec:cde-triangle}

We shall start by recalling the standard setup for the cde-triangle. We let
$(R,\maxid)$ be a local henselian noetherian domain with $K:=\Frac(R)$, the field
of fractions of $R$, and $\k:=R/\maxid$, the residue field of $R$. Recall that
$R$ is \Definition{Henselian} means that for any finite $R$-algebra $A$ any
idempotent of $A/\maxid A$ lifts to an idempotent of $A$. (Note that a field,
which is a local domain with its maximal ideal being the zero ideal, is
Henselian for trivial reasons.) In the case that will eventually interest us $R$
will be the formal power series ring $\pow[\k]t$ in one variable but there are
other choices of interest even in the situation of Kac-Moody algebras. Getting a
little ahead of ourselves we could let $R$ be the algebra of formal power series
expansions around an element of $\h^\vee$, the dual of the Cartan subalgebra of
the Kac-Moody algebra. One could also replace formal power series with
convergent power series if the base field $\k$ is the field of complex
numbers. The reader insisting on having a single example in mind would probably
be best served to let $\k$ be $\C$, the field of complex numbers, and
$R=\pow[\C]{t}$ the ring of formal power series in one variable over $\C$. In
that case $\maxid=tR$. (The case that occurs in the original setup, that of
modular representation theory, is when $R$ equals the ring of $p$-adic integers
or an extension of it obtained by adjoining some roots of unity).

We now assume that we are given an associative $R$-algebra $\cA$ with unit
which is finitely generated and free as an $R$-module. In concrete terms that
means that we are given a basis $e_1=1,e_2,\dots,e_n$ of $\cA$ and structure
constants $a_{ij}^k \in R$ with $e_ie_j=\sum_ka^k_{ij}e_k$ fulfilling the
appropriate conditions making $e_1$ a unit element and makes the multiplication
associative. In the context of representation theory of Lie algebras a good
example to have in mind is the \Definition{Hecke algebra}, $\Cal H$, of a Weyl
group (associated to a finite dimensional semi-simple Lie algebra), where the
parameter $q$ is considered as a formal parameter and $R=\pow[\C]{q-\zeta}$, the
ring of formal expansions around $\zeta$, a root of unity. We get from $\cA$
two algebras over fields, the scalar extension of $\cA$ to $K$, $A:=\Cal
A\Tensor_RK$, and the scalar extension to $\k$, $\overline A:=\Cal
A\Tensor_R\k=A/\maxid A$. Using the concrete description in terms of structure
constants, the first algebra $A$ simply considers the structure constants as
elements of $K$ and for the second algebra one reduces them modulo $\maxid$. We
now assume that $A$ is a semi-simple algebra (i.e., isomorphic to a product of
matrix algebras over some extension field of $K$) but make no such assumption on
$\overline A$. In the case of Hecke algebras it is known that for appropriate
choices of root of unity $\zeta$ one can get examples for which $\overline H$ is
not semi-simple. For the simplest case of a Weyl group of type $A_1$ we have a
basis $1,\sigma$ with $\sigma^2=(q-1)\sigma+q$ and hence for $q=-1$ we get a
non-simple algebra so we can pick $\zeta=-1$.

Recall now that a (finitely generated) $\cA$-, $A$-, or $\overline A$-module $P$ is
\Definition{projective} if for every \Emph{surjective} map \map{g}{M}{N} of modules
and every module homomorphism \map{f}{P}{N} there is a \Emph{lifting}
\map{h}{P}{M} making the diagram
\begin{displaymath}
\xymatrix{
                           &M\ar@{->}[d]^g\\
P\ar@{-->}[ur]^h\ar@{->}[r]_f&N
}
\end{displaymath}
commutative.

The first thing to note is that for the semi-simple algebra $A$ this condition
is trivial; every surjective map \map{g}{M}{N} is \Definition{split}, i.e.,
there is a map of $A$-modules \map{g'}{N}{M} such that $g\circ g'$ is the
identity map. Thus we can find a splitting no matter what $P$ is by putting
$h=g'\circ f$. Hence in the semi-simple case a simple (or irreducible) module
has two properties; it is simple (!) and it is projective. In the non-semi
simple case these two properties are not connected anymore, simple modules are
usually not projective and vice versa. Usually it is the simple modules one is
really interested in but the condition defining projectivity is a very useful
one and the projective modules play at the very least an important auxiliary
r\^ole. More precisely, as we shall see momentarily, the
\Definition{indecomposable} (i.e., those that cannot be written as a non-trivial
direct sum) projective modules are closely related to the simple modules.

An important example of a projective module is the algebra itself. A module
homomorphism from $\cA$ (resp.\ $A$ or $\overline A$) to a module $M$ is
completely determined by where it maps the identity element and any element of
$M$ is such an image. Hence the projectivity follows directly from the
assumption that $M \to N$ is surjective. It is also easy to see that a direct
summand of a projective module is projective. We can write $\cA$ as a direct sum
of indecomposable projective modules and these will then all be projective. We
shall now see that the isomorphism classes of these indecomposable projectives
are in bijection with the isomorphism classes of the simple modules. For this we
need to start by expounding on the significance of the condition that $R$ is
Henselian for indecomposable modules.

Essentially by definition a module $M$ (over $\cA$,\ $A$, or $\overline A$) is
indecomposable precisely when its endomorphism ring contains no non-trivial
idempotents. Now, if $M$ is finitely generated as $R$-module (resp.\ as $K$- or
$\k$-vector space) then so is its endomorphism ring and by the Henselian
property of $R$ its reduction modulo $\maxid$ contains no non-trivial
idempotents. Now, it is a fact that a finite dimensional algebra over a field
with no non-trivial idempotents is \Definition{local}, i.e., if the sum of two
elements is invertible then of the elements is invertible. This property can
then be lifted back to the original endomorphism ring. Note that we only used
that the endomorphism ring was finitely generated so for future use we record
the conclusion as: \Emph{If the endomorphism ring of an indecomposable module is
finitely generated over $R$, $K$, or $\k$, then it is a local ring.} This result
can to begin with then be combined with the Krull-Schmidt theorem which says
that the summands of direct sums of modules whose endomorphism rings are local
can be recovered (up to isomorphism) from the direct sum, i.e., we have
uniqueness of the summands of a direct sum decomposition. This may not be
strictly necessary for our arguments but is certainly reassuring\dots

More seriously we want to use this property to give a bijection between
(isomorphism classes of) indecomposable projective modules and (isomorphism
classes of) simple modules. It follows from the following more precise result.
\begin{proposition}\label{Indecomposable projective}
Let $P$ be a finitely generated indecomposable projective $\cA$-, $A$-, or
$\overline A$-module. 

\part $P$ contains a unique maximal submodule $M_P$ and the
association $P \mapsto P/M_P$ gives a bijection between isomorphism classes of
finitely generated indecomposable projective modules and simple modules.

\part $P$ is isomorphic to a direct summand of $\cA$, $A$, or $\overline A$
respectively.
\begin{proof}
The existence of a maximal submodule follows immediately from the fact that $P$
is finitely generated (as the union of an increasing sequence of proper
submodules can not be all of $P$). Assume therefore that there are two distinct maximal
modules $M,N \subset P$. As they are distinct we get that the ``diagonal'' map
$P \to P/M\Dsum P/N$ is surjective. Using the lifting property of $P$ we may
complete the following diagram with \map{f}{P}{P}
\begin{displaymath}
\xymatrix{
P\ar@{->}[d]\ar@{-->}[rr]^f&&P\ar@{->}[d]\\ 
P/M\Dsum P/N\ar@{->}[r]^-{p_1}&P/M\ar@{->}[r]^-{i_1}&P/M\Dsum P/N,
}
\end{displaymath}
where $p_1$ is the projection on the first factor and $i_1$ is the inclusion in
the first factor. Now, as the endomorphism ring of $P$ is local and as
$1=f+(1-f)$ we have that one of $f$ and $1-f$ is invertible. It is clear that
$1-f$ induces the composite map $P/M\Dsum P/N\mapright{p_2}P/N\mapright{i_2}$
and hence neither $f$ nor $1-f$ can be an isomorphism as if it were, then it
would induce a surjective map $P/M\Dsum P/N \to P/M\Dsum P/N$ which is visibly
not the case.

Assume now that $Q$ is another indecomposable finitely generated projective with
a surjective map $Q \to P/M_P$. By the projectivity of $Q$ there is a lifting of
the map $Q \to P/M_P$ to a map $Q \to P$. This map is surjective as if not its
image is a proper submodule of $P$ and is hence contained in $M_P$ contradicting
the surjectivity of $Q \to P/M/P$. By the projectivity of $P$ the surjective
mapping $Q \to P$ has a splitting $P \to Q$ which makes $P$ a direct summand of
$Q$ which by the indecomposability of $Q$ implies that $P$ is isomorphic to $Q$.

Finally, if $M$ is a simple module, picking $0\ne m \in M$ we get a module
homomorphism $\cA \to M$ (resp\dots) taking $1$ to $m$. Its restriction to
some indecomposable summand $P$ of must be non-zero giving a non-zero map $P \to
M$ but as $M$ is simple this map is surjective. Hence, every simple module $M$
is of the form $P/M_P$ for $P$ a direct summand of $A$. This finishes the proof
of the proposition.
\end{proof}
\end{proposition}
\begin{remark}
Every indecomposable projective module is finitely generated.
\end{remark}
We shall call this indecomposable projective $\cA$-module (resp\dots) which
has a given simple module as quotient the \Definition{projective cover} of the
simple module. In the case of $A$-modules a simple module is its own projective
cover but this is in general not true for $\cA$ or $\oA$. Note now
that every $\oA$-module can be considered as an $\cA$-module by the
(surjective) map $\cA \to \oA$ given by reduction modulo
$\maxid$. This identifies the $\oA$-modules with the $\cA$-modules
that are killed by $\maxid$. Note also that a simple $\oA$-module is
simple as an $\cA$-module. In fact all simple $\cA$-modules $M$ are
obtained in this way. Indeed, $M$ is finitely generated as $\cA$-module as it
is generated by any of its non-zero elements and in particular it is finitely
generated as an $R$-module. Hence Nakayama's lemma says that $\maxid M\ne M$ and
as it is an $\cA$-submodule it must be zero as $M$ is simple. Hence, the
proposition shows that we get a bijection between indecomposable projective
modules over $\cA$ and over $\oA$. We shall see shortly that this
bijection is realised by associating to the projective indecomposable $\Cal
A$-module $P$ the $\oA$-module $\k\Tensor_RP=P/\maxid P$. We start
however by introducing some notation.

A \Definition{$\cA$-lattice} is a $\cA$-module which is finitely generated
and free as an $R$-module. (Hence choosing an $R$-basis of it, the elements of
$\cA$ become represented by matrices with entries in $R$.) Note that
indecomposable projective modules are always lattices as they are direct
summands of $\cA$ and $\cA$ is a lattice (direct summands of finitely
generated free $R$-modules are free as $R$ is local). We then have the following
key result.
\begin{proposition}\label{Homomorphisms}
\part[i] Let $M$ be an $\cA$-lattice and $P$ an indecomposable projective $\Cal
A$-module. Then we have that the $R$-module of $\cA$-homomorphisms
$\Hom_{\cA}(P,M)$ is finitely generated and free. Furthermore, we have that
\begin{displaymath}
K\Tensor_R\Hom_{\cA}(P,M)=\Hom_A(K\Tensor_RP,K\Tensor_RM)
\end{displaymath}
and
\begin{displaymath}
\k\Tensor_R\Hom_{\Cal
A}(P,M)=\Hom_{\oA}(\k\Tensor_RP,\k\Tensor_RM). 
\end{displaymath}
In particular
$\dim_K\Hom_A(K\Tensor_RP,K\Tensor_RM)=\dim_{\k}\Hom_{\oA}(\k\Tensor_RP,\k\Tensor_RM)$.

\part[ii] Assume (for simplicity) that for every simple $A$-module $M$ (resp.\ every
simple $\oA$-module $\overline{M}$) we have that $\End_A(M)=K$ (resp.\
$\End_{\oA}(\overline{M})=\k$). Then for any finitely generated
$A$-module $N$ (resp.\ finitely generated $\oA$-module $\overline{N}$)
and any simple $A$-module (resp.\ simple $\oA$-module) $M$ (resp.\
$\overline{M}$) with projective cover $P$ (resp.\ $\overline{P}$) we have that
$\dim_K\Hom_A(P,N)$ (resp.\
$\dim_\k\Hom_{\oA}(\overline{P},\overline{N})$) is equal to the number
of times that $M$ (resp.\ $\overline{M}$) appears in a Jordan-H\"older sequence
of $N$ (resp.\ $\overline{N}$).
\begin{proof}
We know that $P$ is a direct summand of $\cA$ and the statement behaves well
with respect to taking direct sums. Hence we may replace $P$ by $A$ and then
$\Hom_{\cA}(A,M)=A$ by the map $f \mapsto f(1)$ and similarly for the $A$ and
$\oA$ case. This makes the case $P=A$ of \DHrefpart{i} obvious.

As for \DHrefpart{ii} we prove it by induction over the length of $N$ (the proof
for $\overline N$ is identical) which is finite as $M$ is of finite
dimension. If $N$ is simple then any homomorphism $P \to N$ factors through $P
\to M$ so that $\Hom_A(P,N)=\Hom_A(M,N)$ and by Schur's lemma $\Hom_A(M,N)=0$ if
$M$ and $N$ are non-isomorphic and  $\Hom_A(M,N)$ is $1$-dimensional by
assumption if $M$ and $N$ are isomorphic. This takes care of the case of length
$1$. In the general case choose $N' \subset N$ of length $1$. We then have an
exact sequence
\begin{displaymath}
\shex{\Hom_A(P,N')}{\Hom_A(P,N)}{\Hom_A(P,N/N')},
\end{displaymath}
where the all but the right most exactness is true in general and the right most
exactness is precisely the lifting property for the projective module
$P$. Counting dimensions we get
$\dim_K\Hom_A(P,N)=\dim_K\Hom_A(P,N/N')+\dim_K\Hom_A(P,N')$ which gives the
result by induction.
\end{proof}
\end{proposition}
\begin{remark}
The concrete content of the first part is the following: We may choose an
$R$-basis for $P$ and $M$ and then the action of an element $\cA$ on $P$
resp.\ $M$ are given by matrices with entries in $R$. The corresponding matrices
for $K\Tensor_RP$ and $K\Tensor_RM$ are then obtained by considering the entries
as elements of $K$ and the matrices for $\k\Tensor_RP$ and $\k\Tensor_RM$ are
obtained by reducing the entries modulo $\maxid$. The $R$-module of $\Cal
A$-homomorphisms $P \to M$ then consists of the $R$-matrices that commute with
those of $\cA$ and similarly for the $K$- resp.\ $\k$-vector space of
homomorphism $K\Tensor_RP \to K\Tensor_RM$ resp.\ $\k\Tensor_RP \to \k\Tensor_RM$.

We may choose an $R$-basis for the $R$-module of $\cA$-homomorphisms $P \to
M$ and then the first part says that this basis forms a $K$-basis for the
$K$-vector space of $A$-homomorphisms $K\Tensor_RP \to K\Tensor_RM$ whereas the
reduction modulo $\maxid$ form a basis for the space of
$\oA$-homomorphisms $\k\Tensor_RP \to \k\Tensor_RM$. The $K$-part is
true without assuming that $P$ is projective but the $\k$-part does require that
assumption.
\end{remark}
As a first application we can consider the relation between indecomposable
projective $\cA$-modules and indecomposable projective
$\oA$-modules. Indeed, if $P$ is an indecomposable projective $\Cal
A$-module then $\overline{P}:=P/\maxid P$ is an indecomposable projective
$\oA$-module. That it is projective is clear as any $\Cal
A$-homomorphism $P \to M$ where $M$ is an $\oA$-module factors through
$P \to \overline{P}$ so that the lifting property for $P$ implies that for
$\overline{P}$. On the other hand, it follows from the proposition that
$\End_{\oA}(\overline{P},\overline{P})$ is equal to $\End_{\Cal
A}(P,P)/\maxid \End_{\cA}(P,P)$ and the Henselian property of $R$ implies
that as $\End_{\cA}(P,P)$ has no non-trivial idempotents neither does
$\End_{\oA}(\overline{P},\overline{P})$ so that $\overline{P}$ is
indeed indecomposable.

For the rest of this section we assume that \Emph{the maximal ideal $\maxid$ of
$R$ is generated by a single element} (the standard example of $R=\pow[\k]{t}$
fulfills this condition as does the original example of $p$-adic numbers). What
it means for us is that any finitely generated $R$-submodule of a $K$-vector
space is free which makes it much easier to construct $\cA$-lattices.

Let now $M_1,\dots,M_k$ be a complete set of simple $A$-modules (up to
isomorphism) and $\overline{M}_1,\dots,\overline{M}_\l$ a complete set of simple
$\oA$-modules. We shall now define three integer matrices, the
\Definition{decomposition matrix} $D$, the \Definition{Cartan matrix} $C$, and a
matrix $E$ (which doesn't seem to have acquired a standard name). For a simple
$A$-module (resp.\ $\oA$-module) and an $A$-module (resp.\ $\oA$-module) $N$ of
finite length we put $[N:M]$ equal to the number of times $M$ appears in a
Jordan-H\"older sequence of $N$. We let $D_{ij}$ for $1\le i \le \l$ and $1\le j
\le k$ in the following way. We choose an $\cA$-lattice $\Cal M_j \subset M_j$
such that $M_j =K\Cal M_j$. Such a lattice is easily constructed by taking a
basis of $M_j$ and letting $\Cal M_k$ be the $\cA$-module generated by the
basis. We then let $D_{ij}:=[\k\Tensor_R\Cal M_j:\overline{M}_i]$. By giving a
different formula for this integer we shall see that it is independent of the
choice of $\Cal M_j$ but this can also be shown directly. (Note that the module
$\k\Tensor_R\Cal M_j$ depends in general on the choice of $\Cal M_j$, the
statement is that the components of a Jordan-H\"older is independent of such a
choice.) We define the matrix $E_{ji}$ for $1\le i \le \l$ and $1\le j \le k$ as
follows: We consider the projective cover $P_i$ of $\overline{M}_i$ considered as
an $\cA$-module and put $E_{ji}:=[K\Tensor_RP_i:M_j]$. Finally, We define
$C_{ij}$ for $1\le i,j \le \l$ as $[\overline{P}_j:\overline{M}_i]$. The basic result
concerning these matrices is the following.
\begin{proposition}\label{cde-triangle}
Assume that the endomorphism rings of the $M_j$ and the $\overline{M}_i$
are equal to the base fields ($K$ resp.\ $\k$).

\part We have that $E=D^t$, the transpose matrix, i.e., $[\k\Tensor_R\Cal
M_j:\overline{M}_i]=[K\Tensor_RP_i:M_j]$, where $\Cal M_j$ is an $\cA$-lattice
such that $K\Tensor_R\Cal M_j=M_j$ and $P_i$ a $\cA$-projective cover of $M_i$.

\part We have that $C=DD^t$, i.e., 
\begin{displaymath}
[\overline{P}_j:\overline{M}_i]=\sum_k[\k\Tensor_R\Cal
M_j:\overline{M}_k][K\Tensor_RP_i:M_k].
\end{displaymath}
\begin{proof}
We get from Proposition \ref{Homomorphisms} that 
\begin{displaymath}
[\k\Tensor_R\Cal
M_j:\overline{M}_i]=\dim_\k\Hom_{\oA}(\overline{P}_i,\k\Tensor_R\Cal
M_j)=\dim_K\Hom_A(K\Tensor_RP_i,M_j)
\end{displaymath}
but as $A$ is semi-simple we have that $K\Tensor_RP_i$ is the direct sum of
simple modules and consequently $\dim_K\Hom_A(K\Tensor_RP_i,M_j)$ is equal to
the number of times $M_j$ occurs in $K\Tensor_RP_i$, i.e., is equal to
$[(K\Tensor_RP_i:M_j]$.

As for the second part we have that $[\overline{P}_j:\overline{M}_i]
=\dim_\k\Hom_{\oA}(\overline{P}_i,\overline{P}_j)$ and that is equal to
$\dim_K\Hom_A(K\Tensor_RP_i,K\Tensor_RP_j)$. Writing $K\Tensor_RP_i$ and
$K\Tensor_RP_j$ is a direct sum of $M_k$'s where $M_k$ occurs with multiplicity
$[K\Tensor_RP_i:M_k]$ resp.\ $[K\Tensor_RP_j:M_k]$. This gives that
$\dim_K\Hom_A(K\Tensor_RP_i,K\Tensor_RP_j)$ is equal to
$\sum_k[K\Tensor_RP_i:M_k][K\Tensor_RP_j:M_k]$ which is equal to
$\sum_k[\k\Tensor_R\Cal M_i:\overline{M}_k][\k\Tensor_R\Cal M_j:\overline{M}_k]$
by the first part which proves the result.
\end{proof}
\end{proposition}
\begin{example}
Consider the case of the Hecke algebra of type $A_2$. It has two generators
$\sigma$ and $\tau$ with relations $(\sigma-q)(\sigma+1)=0=(\tau-q)(\tau+1)$ and
$\sigma\tau\sigma=\tau\sigma\tau$. Furthermore, it is $6$-dimensional with a
basis with elements $1$, $\sigma$, $\tau$, $\sigma\tau$, $\tau\sigma$, and
$\sigma\tau\sigma$. (Note also that the Dynkin diagram has a binary symmetry
which gives an automorphism of the Hecke algebra exchanging $\sigma$ and $\tau$
which can be used to simplify the statements to follow.) For any $q$ it has two
$1$-dimensional representations given by letting both $\sigma$ and $\tau$ act by
$-1$ resp.\ by $q$. For general $q$ (more precisely when $q\ne 0,-1,e^{\pm2\pi
i/3}$) this algebra is semi-simple and has these two $1$-dimensional irreducible
representations and one $2$-dimensional. If we let $q \mapsto e^{\pm2\pi i/3}$
then the two $1$-dimensional representations are the \Emph{only} irreducible
representations. Furthermore, their projective covers are $3$-dimensional. Hence
if we let $R=\pow[\C]{q-\zeta}$, where $\zeta = e^{\pm2\pi i/3}$ we have two
projective modules $P_1$ and $P_2$ which are the projective covers of the two
modules $M_1$ and $M_2$ of rank $1$. The kernels of $P_i \to M_i$ are of rank
$2$ isomorphic to the irreducible $2$-dimensional module when scalars are
extended to $K$. That gives the following form for the $E$-matrix
\begin{displaymath}
E=
\left(\!\begin{array}{cc}
1&0\\
0&1\\
1&1
\end{array}\!
\right)
\end{displaymath}
and then the decomposition matrix $D$ is its transpose whereas the Cartan matrix
is
\begin{displaymath}
E^tE=\left(\!
\begin{array}{cc}
2&1\\1&2
\end{array}\!
\right).
\end{displaymath}
If one instead considers the case of $R=\pow[\C]{q}$, then for $q=0$ we have $4$
$1$-dimensional representations for $\sigma,\tau=0,1$. The two representations
where $\sigma,\tau=0$ or $\sigma,\tau=1$ are also projective while the
projective covers of the representations where $\sigma+\tau=1$ have projective
covers of dimension $2$. Note that these two latter projective covers over
$\pow[\k]{q}$ are two distinct lattices (as they are the covers of two distinct
simple modules) yet their $K$-linear extension to representations over $\k((q))$,
the field of formal Laurent series, are isomorphic. This gives an example of two
non-isomorphic lattices in the same representation over $\k((q))$. Their
reductions modulo $q$ are still non-isomorphic but they do indeed have the same
components of their Jordan-H\"older sequences. In any case we have
\begin{displaymath}
E=
\left(\!
\begin{array}{cccc}
1 & 0 & 0 & 0 \\
0 & 1 & 0 & 0 \\
0 & 0 & 1 & 1
\end{array}\!
\right).
\end{displaymath}
\end{example}
\end{section}

\begin{section}{Projective modules over Kac-Moody algebras}

We now intend to do the theory presented in the previous section for modules
over a Kac-Moody algebra. The modules that we shall consider will only very
rarely be finitely generated as $R$-modules. However, going through the previous
section one realises that the important thing is that the $R$-module of
homomorphisms between modules is finitely generated (and at a few points that
the module itself is finitely generated). These two properties will remain true
for the modules that we shall consider. 

Assume that $\k$ has characteristic zero. Let $\g$ be a Kac-Moody algebra over
$R$, which we will take to mean that it is obtained by extending scalars of some
Kac-Moody algebra defined over some subfield of $R$. We put $\g':=\g\Tensor_RK$
and $\overline{\g}:=\g\Tensor_Rk$ and we shall generally let $(-)'$ resp.\
$\overline{(-)}$ denote extension of scalars to $K$ resp.\ $\k$.  By a
\Definition{weight} we shall mean an $R$-linear map from the Cartan subalgebra
of $\g$ to $R$. As usual we have the roots of $\g$ as particular examples of
weights. If $M$ is a $\g$-module and $\lambda$ a weight we define,
unsurprisingly, the weight space
\begin{displaymath}
M_\lambda:= \{m\in M:\forall t\in\h;tm=\lambda(t)m\},
\end{displaymath}
where $\h$ is the Cartan subalgebra. Let $\Gamma$ be a finite set of weights
such that no two distinct elements of $\Gamma+\Delta^ -$ are congruent modulo
$\maxid$, where $\Delta^-$ denotes the set of negative weights.  We then let
$\Cal M_\Gamma$ be the category of finitely generated $\g$-modules $M$ for which
$M$ is the sum of the $M_\lambda$ for $\lambda\in\Gamma+\Delta^-$ and for which
all the $M_\lambda$ are finitely generated $R$-modules. Our assumption on
$\Gamma$ then implies that $M$ is in fact the direct sum of its weight
spaces. We shall say that $M$ is a $\g$-lattice if all the $M_\lambda$ are free
finitely generated $R$-modules. We let $\Gamma'$ resp.  $\overline\Gamma$ denote
the set of induced $\g'$- resp.  $\overline{\g}$-weights and then $\Cal M_{\Gamma'}$
resp.~$\Cal M_{\overline\Gamma}$ have their obvious meaning. We shall now prove
some results imitating \cite{rocha-caridi82::projec+lie}. Recall that a set
$\{P\}$ of projective modules form a \Definition{set of generators} if for every
non-zero module $M$, there is a non-zero homomorphism $P \to M$. In the case of
an algebra $\Cal A$ we have seen that the algebra itself is a generator and that
implies that the set of indecomposable summands of $\Cal A$ form a set of
generators. That they do was a crucial step in showing that every simple module
had a projective cover which was a summand of $\Cal A$.
\begin{lemma}\label{Projective sorites}
\part[i] Let $M,N\in\Cal M_\Gamma$. Then $\Hom_\g(M,N)$ is finitely generated as
$R$-module.

\part[ii] $\Cal M_\Gamma$ has a set of indecomposable projective generators
which all are $\g$-lattices. Reduction modulo $\maxid$ gives a bijection between
the set of isomorphism classes of those and the corresponding set relative to
$\Cal M_{\overline\Gamma}$.

\part[iii] If $P\in\Cal M_\Gamma$ is projective and $M$ is a $\g$-lattice then
$\Hom_\g(M,N)$ is $R$-projective, $\Hom_\g(M,N)' = \Hom_{\g'}(M',N')$ and
$\Hom_\g(M,N)\Tensor_Rk= \Hom_{\overline\g}(\overline M, \overline N)$.
\begin{proof}
As $M$ is finitely generated, any $\g$-map $M\to N$ is determined by its
restriction to a fixed finite set of weight spaces all of which are, by
assumption, finitely generated as R-modules. This observation gives
\DHrefpart{i}. As for \DHrefpart{ii} we follow \cite{rocha-caridi82::projec+lie}
in first constructing a set of projective generators. Let therefore
$\lambda\in\Gamma+\Delta^ -$ and let $N(\lambda)$ be the $\h$-module of rank $1$
corresponding to $\lambda$.  Induce first up to $\b$, the Borel algebra, and
factor out by the submodule consisting of those weight spaces whose weights do
not belong to $\Gamma+\Delta^ -$. This gives a $\b$-module $Q(\lambda)$ finitely
generated free as an $R$-module. Induce then up to $\g$ to get $P(\lambda)$. As
$R$-module $Q(\lambda)$ is a direct summand of $N(\lambda)\Tensor_RU(\n)$ and is
hence $R$-free.  Similarly, $P(\lambda)=Q(\lambda)\Tensor_RU(\n^ -)$ and so is
$R$-free. As $Q(\lambda)$ is a sum of weight spaces, so is $P_\lambda$ and
therefore $P_\lambda$ is a $\g$-lattice being generated by
$Q(\lambda)$. Finally, by construction, for every $\g$-module $M$ in $\Cal
M_\Gamma$ and every $m \in M_\lambda$ there is a unique map of $\h$-modules
$N(\lambda) \to M$ taking the generator $v$ of $N(\lambda)_\lambda$ to $m$. It
induces a map $Q(\lambda) \to M$ of $\b$-modules taking $1\tensor v$ to $m$ and
finally inducing to $\g$ gives a map $P(\lambda) \to M$ taking
$1\tensor(1\tensor v)$ to $m$. Thus a map $P(\lambda) \to M$ is the same thing
as an element $m \in M_\lambda$ and hence $P(\lambda)$ is projective as a
surjective map $M \to N$ induces a surjective map $M_\lambda \to
N_\lambda$. They also form a set of generators as if $M \ne 0$ we have that
$M_\lambda \ne 0$ for some $\lambda$. Being finitely generated $P(\lambda)$ is a
direct sum of indecomposables and as the $P(\lambda)$ form a set of generators
so do the set of such indecomposable summands of some $P(\lambda)$. 

The rest of the proof is almost identical to the proofs of the corresponding
results of Section \ref{sec:cde-triangle}. Hence for instance every indecomposable projective
module contains a maximal submodule as it is finitely generated and the
reduction modulo $\maxid$ of an indecomposable projective is still
indecomposable by the liftability of idempotents as the endomorphism ring is
finitely generated as an $R$-module. The details are left to the reader.
\end{proof}
\end{lemma}
Let us introduce some notation. For each $\k$-weight $\lambda$ we let $\ovl
Z(\lambda)$ denote the corresponding Verma module, $\ovl V(\lambda)$ the
irreducible $\overline{\g}$-module of highest weight $\lambda$ and, supposing
$\lambda\in\overline\Gamma+\Delta^ -$, we let $\ovl I(\lambda)$ the indecomposable
projective, in $\Cal M_{\overline\Gamma}$, which is the projective cover of
$V(\lambda)$. Similarly, the variations $Z(\lambda)$, $I(\lambda)$,
$Z(\lambda)'$ and $V(\lambda)'$ should be self-explanatory.
\end{section}
\begin{section}{The duality theorem}

We are now going to prove the duality theorem. The idea is to consider a
$1$-parameter family of weights which for a general value of the parameter give
a category of modules that is semi-simple, i.e., every module is a direct sum of
simple modules, for which the simple modules are the Verma modules. The Verma
modules for a special value of the parameter then appears as the reduction
modulo $\maxid$ of a Verma module over $R$, which in turn is an $R$-lattice
insider of the simple Verma module for the general value of the weight. This
shows that the decomposition matrix describes the simple components of a Verma
module. On the other hand Rocha-Caridi gives a (more or less) explicit
filtration of an indecomposable projective module whose successive quotients are
Verma modules. This works \Emph{mutatis mutandis} in the $R$-case and then gives
a Jordan-H\"older sequence when scalars are extended to $K$ and hence give the
$E$-matrix. Thus interpreted the duality theorem is simply the statement that $E=D^t$.

We now specialise to $R=\pow[\k]{t}$ and also assume that $\g$ is induced
from $\k$, i.e., that $\g=\ovl\g\Tensor_\k \Pow t$. Let $\ovl \Gamma$ be a finite
set of $\k$-weights. We lift $\ovl\Gamma$ to a set of $R$-weights as
follows. Choose a basis $\{e_i\}$ of $\h$ and $f_i\in t\pow[\k] t$ which are
algebraically independent over $\k$ ($\k((t))$ has infinite transcendence degree
over $\k$ so this is always possible) and let $\Lambda$ be the $R$-weight for
which $\Lambda(e_i)=f_i$. Put now $\Gamma:=\{\gamma+\Lambda:\gamma\in
\overline\Gamma\}$ where we are considering $\k$-weights as $R$-weights in the
obvious way. As the $f_i$ belongs to $t\pow t$ the reduction mod $t$ of $\Gamma$
gives $\ovl\Gamma$ showing the consistency of notation.
\begin{lemma}\label{Generically semi-simple}
The category $\Cal M_{\Gamma'}$ is semi-simple with the Verma modules
$Z(\lambda)'$ for $\lambda\in\Gamma'+\Delta^ -$ as irreducible objects, i.e.,
every module in $\Cal M_{\Gamma'}$ is a direct sum of Verma modules.
\begin{proof}
According to \scite[Thm.\ 1]{kac79::struc+lie} there is a countable set of
hyperplanes in ${\h}^\vee$ such that if $\lambda$ is a $K$-root not lying on any
of them, then the canonical pairing \cite{kac79::struc+lie} on $Z(\lambda)$ is
non-degenerate and as its radical is the maximal proper submodule, $Z(\lambda)$
is irreducible. Now, as any $\lambda\in\Gamma'+\Delta^ -$ has $\k$-algebraically
independent coordinates it therefore follows that $Z(\lambda)'$ is irreducible,
i.e., it equals $V(\lambda)'$. To show semi-simplicity it therefore remains to
show that any extension
\begin{displaymath}
\shex{Z(\lambda)'}{M}{Z(\mu)'}
\end{displaymath}
is trivial. If $\mu \notin \lambda+(\Delta^ -\setminus\{0\})$, then the highest
weight vector of $Z(\mu)'$ lifts to a highest weight vector of $M$ and the
extension splits. If $\mu\in\lambda+(\Delta^ -\setminus\{0\})$ then we take
duals. Here the dual $N^\vee$ of an $N\in\Cal M_{\Gamma'}$ means the following:
Take the set of linear maps $N\to K$ which vanish on all but a finite set of the
weight spaces and consider it as a $\g$-module through the standard action
composed with the canonical involution of $\g$ (which takes a root to its
negative). In this way the dual becomes an involutive anti-equivalence of $\Cal
M_{\Gamma'}$ to itself and the dual of a Verma module is easily seen to be equal
to itself. Therefore, $M^\vee$ is an extension of $Z(\lambda)'$ by $Z(\mu)'$
and as $\lambda \notin \mu+(\Delta^ -\setminus\{0\})$ the extension splits and
by duality so does the original one.
\end{proof}
\end{lemma}
\begin{proposition}\label{Successive Verma quotients}
For any $\lambda\in\Cal M_\Gamma$, $I(\lambda)$ has a finite filtration the
successive quotients of which are Verma modules.
\begin{proof}
The proof of \cite{rocha-caridi82::projec+lie} goes through without
changes. Indeed, using the notation of Lemma \ref{Projective sorites} we can
filter $Q(\lambda)$ with successive quotients being free of rank $1$ as
$R$-modules and with $\n$ acting trivially and inducing such a module to $\g$
gives a Verma module.
\end{proof}
\end{proposition}
Let us now, following \cite{rocha-caridi82::projec+lie}, put $[\ovl
Z(\lambda):\ovl V(\mu)]$ equal to the multiplicity of $\ovl V(\mu)$ in a
composition series for $\ovl Z(\lambda)$ in $\Cal M_{\overline\Gamma} $ and
$[\ovl I(\lambda):\ovl Z(\mu)]$ equal to the multiplicity of $\ovl Z(\mu)$ in
a series as in Lemma \ref{Successive Verma quotients} (with $\Gamma$ replaced by
\ovl\Gamma).
\begin{theorem}
(cf.\ \cite[Thm.\ 4]{rocha-caridi82::projec+lie}) For any
$\lambda,\mu\in\overline\Gamma+\Delta^ -$ we have $[\ovl Z(\Lambda):\ovl
V(\mu)] =  [\ovl I(\mu):\ovl Z(\Lambda)]$.
\begin{proof}
As $\End_{\overline{\g}}(\ovl V(\mu))=\k$ we get that $[\ovl Z(\lambda):\ovl V(\mu)] =
\dim_k\Hom_{\overline{\g}}(\ovl I(\mu),\ovl Z(\lambda))$. By Lemma \ref{Projective
sorites}:iii) this equals $\dim_K\Hom_{\g'}(I(\mu)',Z(\lambda)')$, which by Lemma
\ref{Generically semi-simple} equals the number of times $Z(\lambda)'$ occurs in
$I(\mu)'$. By Lemma \ref{Successive Verma quotients} this equals the number of
times $\ovl Z(\lambda)$ occurs in $\ovl I(\lambda)$, i.e., $[\ovl I(\mu):\ovl
Z(\lambda)]$.
\end{proof}
\end{theorem}
\end{section}
\bibliography{preamble,abbrevs,algebra}

\newcommand\eprint[1]{\texttt{arXiv:#1}}\def\cprime{$'$}
\providecommand{\bysame}{\leavevmode\hbox to3em{\hrulefill}\thinspace}
\providecommand{\MR}{\relax\ifhmode\unskip\space\fi MR }
\providecommand{\MRhref}[2]{%
  \href{http://www.ams.org/mathscinet-getitem?mr=#1}{#2}
}
\providecommand{\href}[2]{#2}
\begin{thebibliography}{RCW82}

\bibitem[KK79]{kac79::struc+lie}
Viktor~G. Kac and David~A. Kazhdan, \emph{Structure of representations with
  highest weight of infinite-dimensional {L}ie algebras}, Adv. in Math.
  \textbf{34} (1979), no.~1, 97--108.

\bibitem[RCW82]{rocha-caridi82::projec+lie}
Alvany Rocha-Caridi and Nolan~R. Wallach, \emph{Projective modules over graded
  {L}ie algebras. {I}}, Math. Z. \textbf{180} (1982), no.~2, 151--177.

\bibitem[Ser78]{serre78::repres}
Jean-Pierre Serre, \emph{Repr{\'e}sentations lin{\'e}aires des groupes finis},
  revised ed., Hermann, Paris, 1978.

\end{thebibliography}
\bibliographystyle{\PRbibstyle}
\end{document}